\documentclass[12pt]{article}
\usepackage[cp1251]{inputenc}
\usepackage[english, russian]{babel}
\usepackage{amsfonts,amssymb,mathrsfs,amscd}
\usepackage{amsthm}
\usepackage{amsmath}

\numberwithin{equation}{section}

\title{Аналитические решения уравнений свертки на выпуклых множествах со смешанной структурой. II}

\author{ С.~Н.~Мелихов, Л.~В.~Ханина}

\date{}

\newtheorem{theorem}[]{Теорема}
\newtheorem{corollary}[]{Следствие}
\newtheorem{lemma}{Лемма}[section]

\theoremstyle{definition}
\newtheorem{definition}[]{Определение}
\newtheorem{remark}[]{Замечание}

\def\CC{\mathbb C}
\def\NN{\mathbb N}

\begin{document}

\maketitle
\thispagestyle{empty}

\begin{abstract}
Доказаны условия существования линейного непрерывного правого
обратного к сюръективному оператору свертки, действующему в пространствах ростков
функций, аналитических на выпуклых подмножествах комплексной плоскости со счетным базисом окрестностей
из выпуклых областей. Они сформулированы в терминах граничного поведения выпуклых
конформных отображений, связанных с указанными множествами.

\medskip

{\bf Ключевые слова: }
уравнение свертки, пространство ростков аналитических функций,
линейный непрерывный правый обратный
\end{abstract}


\section{Введение}\label{s1}

Эта работа является продолжением статьи \cite{MELKHAN1}, в которой
был получен абстрактный критерий существования линейного непрерывного правого
обратного к сюръективному оператору свертки $T_\mu: A(Q+K)\to A(Q)$.
При этом $Q$ -- выпуклое подмножество $\mathbb C$, имеющее счетный базис окрестностей, состоящий из выпуклых областей,
$K\subset\mathbb C$ -- выпуклый компакт, $\mu$ -- линейный непрерывный функционал на $A(K)$.
Для множества $M\subset\mathbb C$ символ $A(M)$ обозначает пространство ростков всех функций, аналитических
на $M$.
Упомянутый критерий \cite[теорема 1]{MELKHAN1} сформулирован в виде  
условия существования специальных семейств субгармонических функций. 
В настоящей работе доказаны аналитические условия
существования линейного непрерывного правого
обратного к сюръективному оператору $T_\mu: A(Q+K)\to A(Q)$.
Они выражены в терминах граничного поведения выпуклых конформных отображений, определяемых $Q$.
Для ограниченного множества $Q$ доказан соответствующий критерий.

Мы придерживаемся обозначений и определений статьи \cite{MELKHAN1}.

\section{Аналитическая реализация абстрактных условий}

\subsection{Характеристики граничного поведения конформных отображений}

Далее нам понадобятся некоторые сведения о выпуклых конформных (биголоморфных)
отображениях и экстремальных субгармонических функциях, связанных с выпуклыми множествами.
Как и ранее, $S:=\{z\in\mathbb C\,|\,|z|=1\}$.

\medskip
\begin{definition} (а) Пусть $G$ -- выпуклая
область в $\CC$, отличная от $\mathbb C$; $\varphi$ -- конформное отображение единичного
круга $\{z\in\CC\,|\,|z|<1\}$ на $G$. Для $r\in(0,1)$ положим
$G_r:=\varphi(\{z\in\CC\,|\,|z|\le r\})$. Из леммы Шварца следует, что
все компакты $G_r$ выпуклы (см. \cite[с.~203]{GOLUSIN}).
Через $H_r$ обозначим опорную функцию $G_r$. Согласно \cite[1.4]{UNIV}
определена функция
\[
D_G(z):=\lim\limits_{r\to
1-0}\frac{H_G(z)-H_r(z)}{1-r}\in(0,+\infty], \,|z|=1.
\]

Пусть $G$ - выпуклый компакт в $\CC$, отличный от точки; $\psi$ -
конформное отображение $\{z\in\CC\,|\, |z|>1\}$ на $\CC\backslash
G$ такое, что $\psi(\infty)=\infty$. Как отмечено в
\cite[замечание 1.3]{SMJ}, для любого $r>1$ компакты
$G_r:=\CC\backslash\psi(\{z\in\CC\,|\,|z|>r\})$ выпуклы. Пусть
$H_r$ - опорная функция $G_r$. Согласно \cite[замечание 4]{MM}
определена функция
\[
D_G(z):=\lim\limits_{r\to
1+0}\frac{H_r(z)-H_G(z)}{r-1}\in[0,+\infty), \,|z|=1.
\]

\medskip
\noindent (б) Следующие определения даны в \cite[\S~1]{ACTA},
\cite[\S\S~1, 2]{CRELLE} (см. также \cite{BORDEAUX}).

Пусть $G\subset\CC$ -- ограниченная выпуклая область,
содержащая $0$. Для $C>0$ через $v_{H_G,C}^{\infty}$ обозначим
наибольшую субгармоническую в $\CC$ функцию такую, что
$v_{H_G,C}^{\infty}\le H_G$ и $v_{H_G,C}^\infty(z)\le C\log|z|
+O(1)$ при $|z|\to \infty$. Отметим, что
$v_{H_G,C}^{\infty}(z)=
Cv_{H_G}^{\infty}(z/C)$
для любых $z\in\CC$, $C>0$. Положим
$v_{H_G}^\infty:=v_{H_G,1}^\infty$. 
Функция $v_{H_G}^\infty$ совпадает с $H_G$ на замкнутом звездном
(относительно $0$) множестве.
Отображение $C_{H_G}^\infty\colon
S\to (0,+\infty]$ определяется равенством
\[
   \left\{z\in\CC\,|\, v_{H_G}^\infty(z)=H_G(z)\right\}=
\left\{\lambda a\,|\,
       a\in S, 0\le \lambda\le 1/C_{H_G}^\infty (a) \right\}.
\]

Пусть теперь $G$ -- выпуклый компакт в $\CC$,
содержащий $0$ в своей относительной внутренности.
Для $c>0$
через $v_{H_G,c}^0$ обозначим наибольшую субгармоническую в
$\CC$ функцию такую, что $v_{H_G,c}^0\le H_G$ и $v_{H_G,c}^0(z)\le
c\,{\rm log}|z|+O(1)$ при $|z|\to 0$. 
Отметим, что
$v_{H_G,c}^{0}(z)=
cv_{H_G}^{0}(z/c)$
для любых $z\in\CC$, $c>0$.
Полагаем
$v_{H_G}^0:=v_{H_G,1}^0$. 
Функция $v_{H_G}^0$ совпадает с $H_G$ на непустом замкнутом
звездном (относительно $\infty$) множестве.
Отображение$C_{H_G}^0\colon S\to[0,+\infty)$
определяется равенством
\[
   \left\{z\in\CC\,|\, v_{H_G}^0(z)=H_G(z)\right\}=\left\{\lambda a\,|\,
       a\in S, 1/C_{H_G}^0(a)\le\lambda<+\infty \right\}.
\]
\end{definition}

Далее вместо $C_{H_G}^0$ и $C_{H_G}^\infty$ будем коротко писать
$C_G^0$ и $C_G^\infty$.

\medskip
Пусть $G$ -- ограниченное выпуклое в $\CC$ множество, внутренность
которого содержит $0$. Согласно \cite[теорема 1.14, теорема 1.18]{ACTA},
\cite[теорема 2.4, теорема 2.13]{CRELLE} (см. также
\cite[теорема 2.6]{BORDEAUX})
найдутся постоянные $b,B>0$, для которых
\[
C_{{\rm int}\,G}^\infty\le D_{{\rm int}\,G}\le BC_{{\rm int}\,G}^\infty, 
\,\,\, bC_{\overline G}^0\le
D_{\overline G}\le C_{\overline G}^0.
\]


Приведем некоторые переформулировки условия $SH(A,Q)$ (см. \cite[п.\,3.3, определение 1]{MELKHAN1})
для специальных случаев.

\begin{remark}
\noindent
 $(i)$ Пусть $Q$ ограниченно, $\omega\ne\emptyset$;\,
$H_n(z):=\max(H_Q(z); H_\omega(z)+|z|/n)$, $z\in\CC$,
$n\in\mathbb N$. Для
любого $A\subset S$ условие $SH(A,Q)$ равносильно следующему:

\noindent
существует семейство
субгармонических в $\CC$ функций $u_t$, $t\in\Gamma(A)$, таких, что $u_t(t)\geqslant
0$, $t\in \Gamma(A)$, и $\forall k \,\exists k'\, \forall s \,\exists s':$
$$
u_t(z)\le H_{k}(z) - H_{k'}(t)-|z|/s'+|t|/s, \,\,
z\in\CC, \, t\in\Gamma(A).
$$
Заметим, что в этом случае $H_n=H_{Q_n}$, $n\in\mathbb N$.



\noindent
$(ii)$ Для $A\subset S$ введем условие $SH_0(A,Q)$
( здесь мы не предполагаем, что $\omega\ne \emptyset$): 

\noindent
Существует семейство субгармонических в $\CC$ функций 
$u_t$, $t\in\Gamma(A)$, таких, что $u_t (t) \ge
0$, $t\in\Gamma(A)$,
и для любого $m \in \NN$ существует $m' \in \NN$, для которого
$$
u_t (z) \le H_{G_{m'}} (z) - H_{G_m} (t), \, z\in\CC, t\in
\Gamma (A).
$$

Пусть $Q$ открыто и отлично от $\CC$.
По \cite[предложение 6, леммы 7, 8]{MOMM94} для любого компакта
$A\subset S$ функция $D_Q$ ограниченна на $A$ тогда и только тогда, 
когда выполняется условие $SH_0(A, Q)$.
\end{remark}

\subsection{Критерий существования специальных семейств
субгармонических функций}

\medskip
Как и ранее, $\omega:=Q\cap\partial Q$, \,
$\omega_0:=(\partial Q)\backslash\omega$.
Введем
множества опорных точек и "биопорных" \, направлений для выпуклого множества $G\subset\CC$,
отличного от $\CC$.
Для $\sigma\subset S$
\[
F_\sigma(G):=\{w\in \partial G\,|\,{\rm Re}(wa)=H_G(a) \mbox{ для
некоторого } a\in\sigma\};
\]
\[
\widehat{\sigma}(G):=\{a\in S\,|\,{\rm Re}(wa)=H_G(a) \mbox{ для
некоторого } w\in F_{\sigma}(G)\}.
\]
Ясно, что $\sigma\subset\widehat{\sigma}(G)$. Если $G=Q$, то будем писать
$F_\sigma$ и $\widehat{\sigma}$ вместо $F_\sigma(G)$ и $\widehat{\sigma}(G)$
соответственно.

Далее $S_0:= S\backslash S_\omega$. Для ограниченного $Q$ полагаем $H_n:=H_{Q_n}$,
$n\in\mathbb N$.

\medskip
\begin{remark} Пусть $Q$ ограниченно.

\noindent
(i) Если
$(\omega_n)_{n\in\NN}$ -- компактное исчерпывание $\omega_0$, то
$(S_{\omega_n})_{n\in\NN}$ -- компактное исчерпывание $S_0$.

\noindent (ii) Если $\kappa$ -- компакт в $S_0$, то
$\widehat{\kappa}$ -- также компакт в $S_0$.

\noindent (iii) Пусть $\kappa, \sigma\subset S$.
Множество $\widehat{\kappa}\cap\sigma$ пусто тогда и только
тогда, когда пусто множество $\kappa\cap\widehat{\sigma}$.

\noindent (iv) Для любого компакта $\kappa$ в $S_0$ существует
$m\in\NN$ такое, что $H_Q=H_m$ на $\kappa$.

\noindent (v) Для любого $k'\in\NN$ существует $k''\in\NN$, для
которого $H_{k''}\le (H_Q+H_{k'})/2$ в $\mathbb C$.

\noindent
(vi) Для любого $k\in\mathbb N$ на $S_\omega$ выполняется строгое неравенство
$H_Q<H_k$.
\end{remark}

\begin{proof}

Заметим, что выпуклое множество $Q_0:=({\rm int}\, Q)\cup\omega_0$
локально замкнуто и строго выпукло в $\partial_r\omega_0$
(см. \cite{Melimo}). Это означает, что множество $Q_0\cap\partial Q_0$
(относительно) открыто в $\partial Q_0$ и пересечение любой
опорной к $\overline{Q_0}$ прямой с множеством $Q_0$ компактно.
Для множеств с такими свойствами утверждения (i) и (ii)
доказаны в \cite{Melimo}.
Поэтому свойства (i), (ii) вытекают из \cite[леммы 3.4, 3.5]{Melimo}.

(iii): Достаточно показать, что из непустоты одного из этих множеств
следует непустота другого. Пусть существует $a\in \widehat{\kappa}\cap\sigma$.
Тогда $a\in\sigma$ и найдется $z\in F_\kappa$ такое, что
${\rm Re}(za)=H_Q(a)$. Вследствие $z\in F_\kappa$ существует
$b\in\kappa$, для которого ${\rm Re}(zb)=H_Q(b)$. Поэтому
$b\in\widehat{\sigma}$, а значит, $b\in\kappa\cap\widehat{\sigma}$.

(iv): Поскольку любая опорная прямая к $\overline Q$ не может
одновременно пересекать и $\omega_0$, и $\omega$,
то $H_Q(a)>H_\omega(a)$ для любого $a\in\kappa$. Значит,
$\alpha:=\inf\limits_{a\in\kappa}(H_Q(a)-H_\omega(a))>0$.
Если $m\in\mathbb N$ таково, что $1/m<\alpha$, то для любого $a\in\kappa$
$$
H_m(a)=\max(H_Q(a); H_\omega(a)+1/m)=H_Q(a).
$$

(v):  Поскольку $H_{k'}=\max(H_Q; H_\omega+1/{k'})$ на $S$,
то можно взять $k'':=2 k'$.

(vi): Пусть $a\in S_\omega$. Тогда
$H_Q(a)=H_\omega(a)$, а значит,
$$
H_k(a)=\max(H_Q(a); H_\omega(a)+1/k)> H_Q(a).
$$

\end{proof}

Докажем далее аналитический критерий для свойства $SH(A,Q)$
для ограниченного $Q$. При этом существенно используются конструкции
З.\,Момма \cite{ACTA}, \cite{CRELLE} (см. определение 1). Метод "склеивания"\, двух семейств субгармонических функций,
примененный при доказательстве импликации (ii)$\Rightarrow$(i) в лемме 2.1,
ранее уже использовался при исследовании операторов свертки \cite{Melimo} и рядов экспонент для аналитических функций
\cite{LOBACH}, \cite{VMJ11}.

\medskip
\begin{lemma} Пусть $Q$ ограниченно, $\omega \neq \varnothing$, 
$A$ -- компактное подмножество $S$. 
Следующие утверждения равносильны:

\begin{itemize}
\item[(i)]
Cуществуют субгармонические в $\mathbb C$ функции $u_t$, $t\in\Gamma(A)$, такие,
что $u_t(t)\ge 0$, $t\in\Gamma(A)$, и $\forall k$ $\exists k'$
$\forall s$ $\exists s'$
$$
u_t(z)\le H_k(z) - H_{k'}(t) -|z|/s' + |t|/s,\,\, z\in\CC,
t\in\Gamma(A).
$$
\item[(ii)]
Функция $D_{\rm int Q}$ ограниченна на каждом компакте в $A\cap
S_0$ и функция $1/D_{\overline Q}$ ограниченна на некоторой
окрестности $A\cap S_\omega$ в $A$.
\end{itemize}
\end{lemma}

\begin{proof} Без ограничения общности можно считать, что
$0\in {\rm int}\,Q$. 

(i) $\Rightarrow$ (ii)): Выберем $k'$ для $k=1$ по (i). По замечанию
2~(vi) на компакте
$S_\omega$ выполняется строгое неравенство $H_Q<H_{k'}$. Поэтому
существуют компактная окрестность $W$ множества $S_\omega$ в $S$ и
$\varepsilon>0$ такие, что на $W$
\begin{equation}
H_{k'}-2\varepsilon \ge H_Q.
\end{equation}
Положим
\[
u:=\left(\sup_{t\in A\cap W}(u_t + H_Q(t))\right)^*,
\]
где $f^*$ обозначает полунепрерывную сверху регуляризацию $f$
(см. \cite[гл.~I, \S~5]{RONKIN}).
Функция $u$ субгармонична в $\CC$, удовлетворяет
неравенству $u\ge H_Q$ на $A\cap W$ и следующему
условию: $\forall n$
$\exists n'$ $\forall m$:
$$
 u(z)\le H_n(z) + \sup\limits_{t\in A\cap W}(H_Q(t)-H_{n'}(t))+1/m, \,\,
z\in \mathbb C.
$$
Следовательно, для любого $n\in\NN$ существует
$n'\in\NN$ такое, что
\begin{equation}
u(z)\le H_n(z) + \sup\limits_{t\in A\cap W}(H_Q(t)-H_{n'}(t)), \,\,
z\in\mathbb C.
\end{equation}
Значит, для любого $n\in\NN$
$$
u(z)\le H_n(z),\,\, z\in\mathbb C.
$$
Переходя к пределу при $n\to\infty$, получим, что $u\le H_Q$ в
$\mathbb C$, а значит, $u=H_Q$ на $A\cap W$. Из оценок (2.1) и (2.2) при
$n=1$ следует, что $u\le -\varepsilon$ в некоторой окрестности
$0$ (тогда $n'$ совпадает с $k'$,
выбранным выше). По лемме 2.14 \cite{CRELLE} для некоторого $\delta_0>0$
выполняется неравенство $v_{H_Q,\delta_0}^0\ge u$ на $A\cap W$.
Поэтому функция $1/C_{\overline Q}^0$ ограниченна на $A\cap W$ (это следует из
ее определения). Значит, и $1/D_{\overline Q}$
ограниченно на $A\cap W$.

\medskip
Пусть множество $\kappa\subset A\cap S_0$ компактно.
По замечанию 2~(ii) множество $\widehat{\kappa}$ -- компакт в $S_0$.
Определим функцию
$$
v:=\left(\sup_{t\in\widehat{\kappa}}(u_t + H_Q(t))\right)^*.
$$
Функция $v$ субгармонична в $\CC$, удовлетворяет неравенству $v\ge
H_Q$ на $\kappa$. Кроме того, выполняется следующее условие: 
$\forall k$ $\exists k'$ $\forall s$
$\exists s'$:
$$
v(z) \le H_k(z) -H_{k'}(t) + H_Q(t) - |z|/s' + 1/s, \,\, z\in\CC, \, t\in A.
$$
Отсюда следует, что $v\le H_Q$ на $\CC$. Поэтому $v=H_Q$ на $\kappa$.

Поскольку $\widehat{\kappa}$ -- компакт в $S_0$, то по замечанию
2\,(iv) $H_Q=H_k$ на $\widehat{\kappa}$ для некоторого $k\in\NN$.
Определим $k'\ge k$ по (i). Пусть $s'$ выбрано в (i) для $s=1$.
Существует окрестность $V$ множества $\widehat{\kappa}$ в
$S$ такая, что
\[
H_k(z)-|z|/s' \le H_Q(z) -|z|/(2 s'),\,\, z\in
\Gamma(V).
\]
Поэтому
\[
v(z) \le H_Q(z) - |z|/(2s') + 1, \,\, z\in\Gamma(V).
\]
Отсюда следует, что $v(z)<H_Q(z)$, если $z\in\Gamma(V)$
и $|z|$ достаточно большое.
Модифицируем $v$ и покажем, что эта модификация удовлетворяет последнему неравенству уже
для всех $z\in\mathbb C$ таких, что $|z|$ достаточно большое.
Положим
\[
L_0(z):=\sup_{w\in F_{\kappa}}{\rm Re}(wz)=\max_{w\in F_{\kappa}}{\rm Re}(wz), \,\, z\in \CC.
\]
Положительно однородная порядка 1 функция $L_0$ удовлетворяет
неравенству $L_0\le H_Q$ на $\CC$ и равенству $L_0=H_Q$ на
$\kappa$. Если $L_0(a)=H_Q(a)$, то существует $w\in F_{\kappa}$,
для которого ${\rm Re}(wa)=H_Q(a)$. Значит, $a\in
S_{F_{\kappa}}=\widehat{\kappa}$. Поэтому $L_0<H_Q$ на
$S\backslash\widehat{\kappa}$. Пусть $\widetilde{v}:=v/2 + L_0/2$.
Субгармоническая в $\CC$ функция $\widetilde v$ удовлетворяет
следующим условиям: \,$\widetilde{v}\le H_Q$ на $\CC$,
$\widetilde{v}=H_Q$
на $\kappa$ и $\widetilde{v}<H_Q$ вне некоторой окрестности начала. По
лемме 2.1 \cite{ACTA} существует $C>0$ такое, что	
$v_{H_Q,C}^\infty\ge \widetilde v$. Поэтому функция $C_{{\rm int}\,Q}^\infty$ (и
$D_{{\rm int}\,Q}$) ограниченна на $\kappa$.

(ii) $\Rightarrow$ (i): Пусть $1/C_{\overline Q}^0$ ограниченно на некоторой
окрестности $\widetilde{A}$ множества $S_\omega\cap A$ в $A$.
Тогда существует $c>0$
такое, что субгармоническая в $\CC$
функция $v_{H_Q,c}^0$ равна $H_Q$
на $\widetilde{A}$. Так как $H_k>H_Q$ на $S_\omega$ для любого
$k\in\NN$, то для каждого $k\in\mathbb N$ найдутся окрестности $S_k$ множества $S_\omega$,
для которых $H_k > H_Q$ на $S_k$.
Поскольку множество $S\backslash S_k$ относительно компактно в $S_0$, то
$\widehat{S\backslash S_k}$ относительно компактно в $S_0$.
Поэтому можно предположить, что $S_{k+1}\subset
S\backslash(\widehat{S\backslash S_k})$ для любого $k\in\NN$ (и
$S_1\cap A\subset\widetilde{A}$). Тогда по замечанию 1~(v) $\widehat{S_{k+1}}\subset S_k$, $k\in\NN$. 
Положим
$A_k:=S_k\cap A$. Поскольку $C_{{\rm int}\,Q}^\infty$ ограниченно на
$A\backslash A_k$ для любого $k\in\NN$, то для каждого $k\in\mathbb N$
существуют постоянные
$C_k>0$ такие, что $v_k^\infty:=v_{H_Q,C_k}^\infty=H_Q$ на
$A\backslash A_{k+2}$. Можно считать, что $C_k\le
C_{k+1}$, $k\in\NN$.

Построим, как и при доказательстве импликации (i)$\Rightarrow$(ii)
(для множества $A_{k+1}$ вместо $\kappa$), положительно однородные
порядка 1 субгармонические в $\CC$ функции $L_k$, не превосходящие
$H_Q$ на $\CC$, равные $H_Q$ на $A_{k+1}$ и такие, что $L_k<H_Q$
на $S\backslash\widehat{A_{k+1}}$. Определим субгармонические
в $\CC$ функции
\[
v_k^0:= v_{H_Q,c}^0/2 + L_k/2, \, k\in\mathbb N.
\]

Зафиксируем $k\in\NN$. 
Заметим, что на $\mathbb C\backslash\{0\}$
\begin{equation}
H_Q<(H_Q-L_k)/2+ H_k.
\end{equation}
Действительно, на $S_k$ выполняется неравенство $H_Q<H_k$. (При этом
$H_Q\ge L_k$ в $\mathbb C$.) На $S\backslash\widehat{A_{k+1}}$ имеет место неравенство
$H_Q>L_k$. Поскольку $\widehat{A_{k+1}}\subset\widehat{S_{k+1}}\subset S_k$, то
$H_Q>L_k$ на $S\backslash S_k\subset S\backslash\widehat{A_{k+1}}$. Это влечет, что
$v_{H_Q,c}^0 <(H_Q-H_k)/2+H_k$ на $\mathbb C\backslash\{0\}$.
Так как $v_{H_Q,c}^0(z)\le c{\rm log}|z|+O(1)$, $|z|\to 0$,
то найдется $\widetilde k\in\mathbb N$, для которого на $\mathbb C$
\[
  v_k^0\le 
L_k/2+(H_Q-L_k)/4+H_k/2-1/\widetilde k=
\]
\[
H_k -D_k/2 -1/\widetilde{k},
\]
где $D_k:= H_k - (H_Q+L_k)/2$.
Вследствие (2.3) в $\mathbb C\backslash\{0\}$ выполняется неравенство $D_k>0$.

Выберем $k'$ так, чтобы $H_{k'}-H_Q\le 1/\tilde{k}$ на $A_{k+1}$.
Тогда для любого $s\in\NN$ найдется $s'\in\NN$, для которого
$D_k(z)/2\ge |z|/s'$, $z\in\CC$, и
\begin{equation}
   v_k^0(z)\le H_k(z) -|z|/s' + H_Q(t) - H_{k'}(t)+1/s,\,z\in\CC,\,
t\in A_{k+1}.
\end{equation}
При этом можно выбрать $k'$ таким большим, что $H_Q=H_{k'}$ на
$S\backslash A_{k+2}$. Поскольку $v_k^\infty(z)\le C_k{\rm
log}|z|+O(1)$, $|z|\to\infty$, и $v_k^\infty\le H_Q$, то для
любого $s$ можно выбрать $s'$ так, чтобы выполнялось и неравенство
$v_k^\infty(z)\le
H_Q(z)-|z|/s' +1/s$, $z\in\CC$. Поэтому
\begin{equation}
   v_k^\infty(z) \le H_k(z)-|z|/s' +H_Q(t)-H_{k'}(t)+1/s,
        \,\, z\in\CC,\, t\in S\backslash A_{k+2}.
\end{equation}
Итак, для любого $k\in\mathbb N$ существует $k'\in\mathbb N$ такое, что
для любого $s\in\mathbb N$ найдется $s'\in\mathbb N$, для которого выполняются неравенства
(2.4) и (2.5).

Заметим, что $v_l^0\ge v_{l+1}^0$ и
$v_{l}^\infty\le v_{l+1}^\infty$ на $\mathbb C$ для любого $l\in\mathbb N$. Значит,
существует (поточечный)
$\lim\limits_{l\to\infty}
v_l^0=:v_\infty^0$. При этом $v_\infty^0$ -- субгармоническая в
$\CC$ функция такая, что $v_\infty^0=H_Q$ на
$A_0:=\bigcap\limits_{l\in\NN}A_l$.

Для $t\in A\backslash A_2$ положим $\tilde{v}_t:=v_1^\infty$. Для
$l\in\NN$, $t\in A_{l+1}\backslash A_{l+2}$ пусть
$\tilde{v}_t:=v_l^0/2+v_l^\infty/2$. Для $t\in A_0$ определим
$\tilde{v}_t:=v_\infty^0$. По построению $\tilde{v}_t(t)=H_Q(t)$
для каждого $t\in S$.

Пусть $t\in A_{l+1}\backslash A_{l+2}$. Для $k\le l$ и $k'$, $s$ и
$s'$, как в (2.4),
\[
   \tilde{v}_t(z)\le
      (H_k(z) -|z|/s' +H_Q(t)-H_{k'}(t)+1/s)/2 + H_Q(z)/2,\, z\in\CC.
\]
Последние неравенства выполняются также для всех $k\in\mathbb N$ 
для функций 
$\tilde v_t:=v_\infty^0$, если $t\in A_0$.
Если $k\ge l$ и $k'$, $s$ и $s'$, как в (2.5), то выполняется
неравенство
\[
    \tilde{v}_t(z)\le
       H_Q(z)/2+(H_k(z)-|z|/s' +H_Q(t)-H_{k'}(t)+1/s)/2,\, z\in\CC.
\]
Эти неравенства имеют место и для функций $\tilde
v_t:=v_1^\infty$, если $t\in A\backslash A_2$.

По замечанию 2\,(v) для любого $k'$ существует $k''$ такое, что
$(H_Q+H_{k'})/2\ge H_{k''}$, а значит, $(H_Q-H_{k'})/2\le
H_Q-H_{k''}$ в $\mathbb C$. Поэтому для всех $t\in A$
\[
    \tilde{v}_t(z)\le
      H_k(z) -|z|/(2s') +H_Q(t)-H_{k''}(t)+1/(2s),\,
\,z\in\CC.
\]
Искомыми являются функции
\[
u_t(z):=|t|\tilde v_{t/|t|}(z/|t|)-H_Q(t),\, z\in\CC,\,
t\in\CC\backslash\{0\}; \,\, u_0\equiv 0.
\]
\end{proof}

\subsection{Основной результат}

По-прежнему $S_0:=S\backslash S_\omega$.

\medskip
\begin{theorem}
	Пусть $Q$ -- выпуклое подмножество $\CC$  с непустой внутренностью,
 обладающее счетным базисом
	окрестностей, состоящим из выпуклых областей; $\omega = Q\cap \partial Q$ непусто;
	$K$ -- выпуклый компакт в $\mathbb C$.
	
    (I) Пусть множество нулей  $\widehat\mu$
    конечно или пусто. Тогда оператор
    $T_\mu: A(Q+K)\to A(Q)$ имеет ЛНПО.

\noindent
(II)  Пусть множество нулей ненулевой функции $\widehat\mu$
бесконечно и оператор
$T_\mu: A(Q+K)\to A(Q)$
сюръективен. Если $Q$ ограниченно, то следующие утверждения равносильны:
\begin{itemize}
\item[(i)] $T_\mu: A(Q+K)\to A(Q)$ имеет ЛНПО.
\item[(ii)]  Функция $D_{\rm int Q}$ ограниченна на каждом компакте в
$A_{\widehat{\mu}}\cap S_0$, а функция $1/D_{\overline Q}$ 
ограниченна на некоторой окрестности $A_{\widehat{\mu}}\cap S_\omega$ в $A_{\widehat{\mu}}$.
\end{itemize}

Пусть $Q$ неограниченно. Зафиксируем произвольную прямую $l$, разделяющую $Q$. 
Утверждение (i) вытекает из такого:

\begin{itemize}
\item[(iii)] Функция $D_{\rm int Q}$ ограниченна на каждом компакте в
$A_{\widehat{\mu}}\cap S_0$, а функция $1/D_{\overline{Q(l)}}$ ограниченна на некоторой
окрестности $A_{\widehat{\mu}}\cap S_\omega$ в $A_{\widehat{\mu}}$.
\end{itemize}

\end{theorem}

\begin{proof}
(I): Найдутся $\lambda\in\CC$ и многочлен $P$ такие, что 
$\widehat{\mu}(z)=e^{\lambda z}P(z)$, $z\in\CC$. Так как оператор 
$T_\mu: A(Q_n+K)\to A(Q_n)$ сюръективен для любого $n\in\NN$, то
$T_\mu$ отображает $A(Q+K)$ на $A(Q)$. По теореме об открытом отображении
для (LF)-пространств отображение
$T_\mu: A(Q+K)\to A(Q)$ открыто. Это и конечномерность ядра оператора $T_\mu$
в $A(Q+K)$ влекут существование ЛНПО к $T_\mu$.

(II): $(i)\Rightarrow(ii)$: 
По \cite[теорема 1]{MELKHAN1} выполняется условие $SH(A_{\widehat\mu}, Q)$.
По лемме 2.1 (с учетом замечания 1~(i))
справедливо утверждение (ii).

Докажем теперь справедливость импликаций (ii)$\Rightarrow$(i) для ограниченного
$Q$ и (iii)$\Rightarrow$(i) для неограниченного $Q$. 
Прежде всего покажем, что в обоих случаях выполнено условие $SH(A_{{\widehat\mu}}, Q)$.

Пусть $Q$ ограниченно. Тогда условие
$SH(A_{{\widehat\mu}}, Q)$ имеет место вследствие леммы 2.1. 

Пусть теперь
$Q$ неограниченно. То, что любая опорная прямая к $\overline Q$ не может
одновременно пересекать и $\omega$, и $(\partial Q)\backslash\omega$ 
\cite[лемма 2.1]{MELKHAN1},
влечет, что
существует окрестность $V$ множества
$S_\omega$ в $S$ такая, что $F_V=F_V(Q(l))$.
По условию (ii) найдется компактная окрестность $W\subset V$ множества
$S_\omega$ в $S$, для которой 
функция $1/D_{\overline{Q(l)}}$ ограниченна на $A_{\widehat\mu}\cap W$.
Пусть $\sigma$ -- компактное подмножество $(A_{\widehat\mu}\cap W)\cap S_0$.
Тогда $D_{{\rm int}\,Q}$ ограниченно на $\sigma$. Вследствие
локальности свойства ограниченности $D_{{\rm int}\,Q}$ \cite[предложение 4.1]
{UNIV} функция $D_{{\rm int}\,Q(l)}$ также ограниченна на $\sigma$.
По лемме 2.1 выполняется условие $SH(A_{\widehat\mu}\cap W, Q(l))$.
Поскольку $H_Q\left|\right._{A_{\widehat\mu}\cap V}= 
H_{Q(l)}\left|\right._{A_{\widehat\mu}\cap V}$ и $H_{Q(l)}\le H_Q$ в $\CC$,
то выполняется также условие $SH(A_{\widehat\mu}\cap W, Q)$.
Пусть ${\rm int}_r W$ обозначает относительную внутренность $W$ в $S$.
Тогда $\tau:=A_{\widehat\mu}\backslash ({\rm int}_r W)$ --
компакт в $S\backslash S_\omega$, а значит, функция
$D_{{\rm int}\,Q}$ ограниченна на $\tau$. По замечанию 1
имеет место условие $SH_0(\tau, {\rm int}\,Q)$. Поэтому справедливо и
условие $SH(\tau, Q)$. Поскольку $A_{\widehat\mu}=(A_{\widehat\mu}\cap W)\cup\tau$,
то имеет место $SH(A_{\widehat\mu}, Q)$.
Согласно \cite[лемма 3.7]{MELKHAN1}   $A_{\widehat\mu}\subset R_Q$. Это влечет, что
функция $\hat\mu$ медленно убывает на
$\Gamma(A_{\hat\mu})$ \cite[лемма 3.4]{MELKHAN1}.
Поэтому существует последовательность областей
$(\Omega_j)_{j\in\mathbb N}$, как в \cite[лемма 3.3]{MELKHAN1}

Далее доказательство наличия ЛНПО к $T_\mu: A(Q+K)\to A(Q)$
стандартно 
(используются построения, предложенные в \cite{MOMM94}; см.
доказательство импликации (ii)$\Rightarrow$(i) в теореме 4.2 в
\cite{BM}).
По  \cite[ лемма 3.6 ]{MELKHAN1} отображение
$$
\rho : A_{Q + K}\left/\right.(\widehat\mu \cdot A_Q) \to \Lambda
(Q, \widehat\mu, \mathbb E), \,\, f + \widehat\mu \cdot A_Q
\mapsto \left(f\left|\right._{\Omega_j}+I_j)\right)_{j\in\mathbb N}
$$
является топологическим изоморфизмом "на". Пусть $e_{jk}$, 
$1\le k\le k_j:={\rm dim}\,E_j$, --
базис Ауэрбаха в $E_j$; $X_{jk}$ -- последовательность
$(0,...,0, e_{jk}, 0,...)$, в которой единственная ненулевая координата
$e_{jk}$ находится на $\left(k+\sum\limits_{i=0}^{j-1} k_i\right)$-ой
позиции ($k_0:=0$). Семейство (последовательность) $\left(X_{jk}\right)_{1\le k\le k_j,
j\in\mathbb N}$ -- абсолютный базис в $\Lambda(Q, \widehat\mu, \mathbb E)$.
Найдутся целые функции $g_{jk}$ такие, что
$$
\left(g_{jk}\left|_{\Omega_l}+I_l\right.\right)_{l\in\mathbb N}= X_{jk},
\,1\le k\le k_j, \, j\in\mathbb N.
$$
Кроме того, $g_{jk}$ удовлетворяют следующим оценкам сверху:
$\forall n$ $\exists n'$ $\forall m$ $\exists m'$ $\exists C_1<+\infty$:
$$
|g_{jk}(z)|\le C_1\exp(H_{n m'}(z) -H_{n' m}(\nu_j)),
\,z\in\mathbb C, \,1\le k\le k_j, \, j\in\mathbb N.
$$
Отсюда следует, что отображение
$$
R\left(\sum\limits_{1\le k\le k_j, j\in\mathbb N} c_{jk} X_{jk}\right):=
\sum\limits_{1\le k\le k-j, j\in\mathbb N} c_{jk}g_{jk}
$$
непрерывно (и линейно) действует из $\Lambda(Q, \widehat\mu, \mathbb E)$
в $A_{Q+K}$ (ряд в определении $R$ абсолютно сходится в $A_{Q+K}$). 
При этом оператор $\rho$ является ЛНПО к фактор-отображению
$q: A_{Q + K} \to A_{Q + K} / (\widehat\mu \cdot
A_Q)$. Кроме того,
$\widehat\mu \cdot A_Q$ замкнуто в $A_{Q+K}$.
Вследствие   \cite[лемма 2.5]{MELKHAN1} 
$T_\mu: A(Q+K)\to A(Q)$ также имеет ЛНПО.
\end{proof}

\medskip
Поскольку существует сюръективный оператор $T_\mu: A(Q+K)\to A(Q)$ 
такой, что $A_{\widehat\mu}=S$, справедливо

\begin{corollary}
Для множеств $Q$ и $K$, как в теореме 1, следующие утверждения равносильны:
\begin{itemize}
\item[(i)] Любой сюръективный оператор
$T_\mu: A(Q+K)\to A(Q)$ имеет ЛНПО.
\item[(ii)] Множество $Q$ ограниченно, функция $D_{\rm int Q}$ ограниченна на каждом компакте в
$S_0$ и функция $1/D_{\overline Q}$ ограниченна на некоторой
окрестности $S_\omega$ в $S$.
\end{itemize}
\end{corollary}

\begin{remark} (i) Пусть $Q$ ограниченно;
$\varphi$ -- конформное отображение единичного
круга $\{z\in\CC\,|\,|z|<1\}$ на ${\rm int}\,Q$, а $\psi$ -- 
конформное отображение $\{z\in\CC\,|\, |z|>1\}$ на $\CC\backslash
\overline Q$ такое, что $\psi(\infty)=\infty$. Согласно 
\cite{UNIV}, \cite{SMJ}
ограниченность $D_{\rm int Q}$ на $S$ равносильна ограниченности
$\varphi'$ на $D_{{\rm int}\, Q}$, а ограниченность $1/D_{\overline Q}$
на $S$ равносильна ограниченности $1/\psi'$ на
$\{z\in\CC\,|\, |z|>1\}$.

\noindent
(ii) Условия ограниченности $D_{\rm int Q}$ и $1/D_{\overline Q}$
(для ограниченного $Q$) в следствии 1 имеют место, если,
например, $\partial Q\in C^{1,\lambda}$ для некоторого $\lambda>0$
(см. \cite[пример 2.8]{BORDEAUX}).
Одно из них не выполняется, если $\partial Q$
имеет угловую точку (см. \cite{UNIV}, \cite{SMJ}).

\end{remark}

\medskip

\begin{flushleft}
С.~Н.~Мелихов

Южный федеральный университет, Ростов-на-Дону;
Южный математический институт Владикавказского научного центра РАН,
Владикавказ

E-mail: melih@math.rsu.ru

\end{flushleft}

\begin{flushleft}
Л.~В.~Ханина

Южный федеральный университет, Ростов-на-Дону

E-mail: khanina.lv@mail.ru

\end{flushleft}

\end{document}